\documentclass[11pt,english]{article}

\usepackage[T1]{fontenc}
\usepackage[latin9]{inputenc}
\usepackage{geometry}
\usepackage{bm}
\usepackage{amsthm}
\usepackage{amsmath}
\usepackage{amssymb}
\usepackage{graphicx}

\makeatletter

\providecommand{\tabularnewline}{\\}

\numberwithin{equation}{section}
\numberwithin{figure}{section}
\theoremstyle{plain}
\newtheorem{thm}{\protect\theoremname}[section]
  \theoremstyle{plain}
  \newtheorem{cor}[thm]{\protect\corollaryname}
  \theoremstyle{definition}
  \newtheorem{defn}[thm]{\protect\definitionname}
  \theoremstyle{plain}
  \newtheorem{lem}[thm]{\protect\lemmaname}
  \theoremstyle{remark}
  \newtheorem{note}[thm]{\protect\notename}

\usepackage[OT1]{fontenc}

\def\trace{\mathop{\operator@font Tr}\nolimits}
\def\diag{\mathop{\operator@font diag}\nolimits}
\def\var{\mathop{\operator@font Var}\nolimits}
\def\ops{\mathop{\operator@font ops}\nolimits}
\def\E{\mathop{\operator@font E}\nolimits}
\def\rank{\mathop{\operator@font rank}\nolimits}
\def\erfc{\mathop{\operator@font erfc}\nolimits}
\def\round{\mathop{\operator@font round}\nolimits}
\def\poly{\mathop{\operator@font poly}\nolimits}
\def\Null{\mathop{\operator@font null}\nolimits}
\def\enull{\mathop{\operator@font enull}\nolimits}
\def\range{\mathop{\operator@font range}\nolimits}
\def\span{\mathop{\operator@font span}\nolimits}

\makeatother

\usepackage{babel}
  \providecommand{\corollaryname}{Corollary}
  \providecommand{\definitionname}{Definition}
  \providecommand{\lemmaname}{Lemma}
  \providecommand{\notename}{Note}
\providecommand{\theoremname}{Theorem}

\begin{document}

\title{\textbf{Effective Stiffness: Generalizing Effective Resistance Sampling
to Finite Element Matrices}}

\author{\textbf{Haim Avron}\\
IBM T.J. Watson Research Center \and \textbf{Sivan Toledo}\\
Tel Aviv University}
\maketitle
\begin{abstract}
We define the notion of \emph{effective stiffness} and show that it
can used to build \emph{sparsifiers}, algorithms that sparsify linear
systems arising from finite-element discretizations of PDEs. In particular,
we show that sampling $O(n\log n)$ elements according to probabilities
derived from effective stiffnesses yields a high quality preconditioner
that can be used to solve the linear system in a small number of iterations.
Effective stiffness generalizes the notion of effective resistance,
a key ingredient of recent progress in developing nearly linear symmetric
diagonally dominant (SDD) linear solvers. Solving finite elements
problems is of considerably more interest than the solution of SDD
linear systems, since the finite element method is frequently used
to numerically solve PDEs arising in scientific and engineering applications.
Unlike SDD systems, which are relatively easy to solve, there has
been limited success in designing fast solvers for finite element
systems, and previous algorithms usually target discretization of
limited class of PDEs like scalar elliptic or 2D trusses. Our sparsifier
is general; it applies to a wide range of finite-element discretizations.
A sparsifier does not constitute a complete linear solver. To construct
a solver, one needs additional components (e.g., an efficient elimination
or multilevel scheme for the sparsified system). Still, sparsifiers
have been a critical tools in efficient SDD solvers, and we believe
that our sparsifier will become a key ingredient in future fast finite-element
solvers. 
\end{abstract}

\section{Introduction}

We explore the sparsification of finite element matrices using \emph{effective
stiffness} \emph{sampling}. The goal of the sparsification is to reduce
the number of elements in the matrix so that it can be easily factored
and used as a preconditioner for an iterative linear solver. We show
that sampling non-uniformly $O(n\log n)$ elements produces a matrix
that is with high probability spectrally close to the original matrix,
and therefore an excellent preconditioner. The sampling probability
of an element is given by the largest generalized eigenvalue of the
element matrix and the effective stiffness matrix of the element. 

Effective stiffness generalizes the notion of effective resistance,
a key ingredient in much of the recent progress in nearly optimal
symmetric diagonally dominant (SDD) linear solvers~\cite{KoutisMillerPang10,BlellochEtAl11,KoutisMillerPang11}.
Solving finite elements problems is of considerably more interest
than the solution of SDD linear systems, since the finite element
method is frequently used to numerically solve PDEs arising in scientific
and engineering applications.

Unlike SDD systems, which are relatively easy to precondition, there
has been limited success in designing fast solvers for finite element
systems. Efforts to generalize combinatorial preconditioners to matrices
that are not weighted Laplacians followed several paths, and started
long before recent progresses. Gremban showed how to transform a
linear system whose coefficient matrix is a signed Laplacian to a
linear system of twice the size whose matrix is a weighted Laplacian.
The coefficient matrix is a $2$-by-$2$ block matrix with diagonal
blocks with the same sparsity pattern as the original matrix $A$
and with identity off-diagonal blocks. A different approach is to
extend Vaidya's construction to signed graphs~\cite{BCHT04}. The
class of symmetric matrices with a symmetric factorization $A=UU^{T}$
where columns of $U$ have at most $2$ nonzeros contains not only
signed graphs, but also gain graphs, which are not diagonally dominant~\cite{BCPT05};
it turns out that these matrices can be scaled to diagonal dominance,
which allows graph preconditioners to be applied to them~\cite{DS08}.

The matrices that arise in finite-element discretization of elliptic
partial differential equations (PDEs) are positive semi-definite,
but in general they are not diagonally dominant. However, when the
PDE is scalar (e.g., describes a problem in electrostatics), the matrices
can sometimes be approximated by diagonally dominant matrices. In
this scheme, the coefficient matrix $A$ is first approximated by
a diagonally-dominant matrix $D$, and then $G_{D}$ is used to construct
the graph $G_{B}$ of the preconditioner $B$. For large matrices
of this class, the first step is expensive, but because finite-element
matrices have a natural representation as a sum of very sparse matrices,
the diagonally-dominant approximation can be constructed for each
term in the sum separately. There are at least three ways to construct
these approximations: during the finite-element discretization process~\cite{BHV08},
algebraically~\cite{ACST09}, and geometrically~\cite{WS06}. A
slightly modified construction that can accommodate terms that do
not have a close diagonally-dominant approximation works well in practice~\cite{ACST09}.

Another approach for constructing combinatorial preconditioners to
finite element problems is to rely on a graph that describes the relations
between neighboring elements. This graph is the dual of the finite-element
mesh; elements in the mesh are the vertices of the graph. Once the
graph is constructed, it can be sparsified much like subset preconditioners.
This approach, which is applicable to vector problems like linear
elasticity, was proposed in~\cite{ShklarskiToledo08}; this paper
also showed how to construct the dual graph algebraically and how
to construct the finite-element problem that corresponds to the sparsified
dual graph. The first effective preconditioner of this class was proposed
in~\cite{DS07}. It is not yet known how to weigh the edges of the
dual graph effectively, which limits the applicability of this method.
However, in applications where there is no need to weigh the edges,
the method is effective~\cite{ShklarskiToledo09}.

Our theory of effective stiffness sampling is an extension of the
theory of effective resistance sampling. It is applicable to a wide
range of finite element discretizations. But our sparsifier is not
yet a complete algorithm for solving finite-element systems. We discuss
the remaining challenges in Section~\ref{sec:discussion}. Nevertheless,
we our results constitute a useful technique that should lead to fast
finite-element solvers. A similar evolution gave rise to the fastest
SDD solvers: Spielman and Srivastava's theory of effective resistance
sampling~\cite{SpielmanSrivastava08} did not immediately lead to
efficient algorithm, but the follow-up work of Koutis et al. turned
it into very efficient algorithms~\cite{KoutisMillerPang10,KoutisMillerPang11}.
The techniques used by the authors of~\cite{KoutisMillerPang10,KoutisMillerPang11}
to solve SDD systems do not trivially carry over to finite element
matrices. For example, their constructions rely on low-stretch trees,
a concept that does not have a natural extension for finite element
matrices. But we expect such extensions to be developed in the future.

\section{Preliminaries}

\subsection{Notation}

We use $[n]$ to denote the set $\{1,\dots,n\}$. We use $A,B,\ldots$
to denote matrices; $x,y,\ldots$ to denote column vectors. $e_{i}$
is the $i$th standard basis vector (whose dimensionality will be
clear from the context, or explicitly stated): all entries all zero
except the $i$th entry which equals one. We denote by $A^{+}$ the
Moore-Penrose pseudo-inverse of $A$. For a symmetric positive definite
matrix $A$, $\lambda_{\max}(A)$ is the maximum eigenvalue, $\lambda_{\min}(A)$
is the minimum eigenvalue and $\kappa(A)$ is the condition number,
that is $\lambda_{\max}(A)/\lambda_{\min}(A)$. For two symmetric
matrices $A$ and $B$ of the same dimension, we denote by $A\preceq B$
that $x^{T}Ax\leq b^{T}Bx$ for all $x$. We abbreviate ``independent
identically distributed'' to ``i.i.d'', ``with probability'' to ``w.p''
and ``with high probability'' to ``w.h.p''.

\subsection{Sums of Random Matrices}

Approximating a matrix using random sampling can be viewed as a particular
case of taking sums of random matrices. In the last few years there
has been significant literature on showing concentration bounds on
such sums~\cite{RudelsonVershynin07,MagenZouzias10,Oliveira10,Tropp11}.
We use the following Matrix Chernoff bound due to Tropp~\cite{Tropp11}.

\begin{thm}
\label{thm:Chernoff}\cite[Theorem~1.1]{Tropp11}Let $A_{1},A_{2},\dots,A_{M}$
be independent matrix-valued random variables. Assume that the $A_{i}$s
are real, $n$-by-$n$ and symmetric positive semidefinite with $\Vert A_{i}\Vert_{2}\leq\gamma$
almost surely for all $i$. Define 
\[
\mu_{\min}=\lambda_{\min}\left(\sum_{i=1}^{M}E(A_{i})\right)\,\,\, and\,\,\,\mu_{\max}=\lambda_{\max}\left(\sum_{i=1}^{M}E(A_{i})\right)\,.
\]
Then for $\eta\in[0,1]$ we have 
\[
\Pr\left(\lambda_{\min}\left(\sum_{i=1}^{M}A_{i}\right)\leq(1-\eta)\mu_{\min}\right)\leq n\left[\frac{\exp(-\eta)}{(1-\eta)^{(1-\eta)}}\right]^{\mu_{\min}/\gamma}
\]
and 
\[
\Pr\left(\lambda_{\max}\left(\sum_{i=1}^{M}A_{i}\right)\geq(1+\eta)\mu_{\max}\right)\leq n\left[\frac{\exp(\eta)}{(1+\eta)^{(1+\eta)}}\right]^{\mu_{\max}/\gamma}\,.
\]

\end{thm}
The following is an immediate corollary.
\begin{cor}
\label{cor:kappa-bound}Let $A_{1},A_{2},\dots,A_{M}$ be independent
matrix-valued random variables. Assume that the $A_{i}$s are real,
$n$-by-$n$ and symmetric positive definite with $\E(A_{i})=I_{n}$
and $\Vert A_{i}\Vert_{2}\leq\gamma$. Let $\kappa_{\max}>1$ and
$\delta\in(0,1)$, and define 
\begin{equation}
C(\kappa_{\max})=\frac{\kappa_{\max}+1}{2\kappa_{\max}\ln(2\kappa_{\max}/(\kappa_{\max}+1))-\kappa_{\max}+1}\,.\label{eq:Ceq}
\end{equation}
If $M\geq C(\kappa_{\max})\gamma\ln(2n/\delta)$ then 
\[
\Pr\left(\frac{1}{M}\sum_{i=1}^{M}A_{i}\,\mbox{\textrm{is singular}}\,\,\mathrm{or}\,\,\kappa\left(\frac{1}{M}\sum_{i=1}^{M}A_{i}\right)>\kappa_{\max}\right)\leq\delta\,.
\]
\end{cor}
\begin{proof}
Use Theorem~\ref{thm:Chernoff} with $\eta=(\kappa_{\max}-1)/(\kappa_{\max}+1)$
to show that all eigenvalues of $\frac{1}{M}\sum_{i=1}^{M}A_{i}$
are smaller than $1-\eta$ with probability at most $\delta/2$ and
bigger than $1+\eta$ with probability of at most $\delta/2$ each.
Union-bound ensures that all eigenvalue are within $[1-\eta,1+\eta]$
with probability of at least $1-\delta$. This establishes the bound
on the condition number with high probability.
\end{proof}

\subsection{Generalized eigenvalues, analysis of iterative methods and sparsification
bounds}

A well known property of many iterative linear solvers, including
the popular conjugate gradient and the theoretically convenient Chebyshev
iteration, is that their convergence rate depends on the distribution
of the eigenvalues of the coefficient matrix (its spectrum). The rate
depends on how much the spectrum is clustered, but it is hard to form
a concise bound. A simple and useful theoretical bound for symmetric
positive semidefinite matrices depends only on the ratio between the
largest and smallest eigenvalue. When using preconditioned methods
convergence is governed by the generalized eigenvalues.
\begin{defn}
\label{def:gev}Given two matrices $A$ and $B$ with the same null
space $\bm{\mathrm{N}}$, a \emph{finite generalized eigenvalue} $\lambda$
of $(A,B)$ is a scalar satisfying $Ax=\lambda Bx$ for some $x\not\in\bm{\mathrm{N}}$.
The \emph{generalized finite spectrum} $\Lambda(A,B)$ is the set
of finite generalized eigenvalues of $(A,B)$. If both $A$ and $B$
are symmetric positive definite,the \emph{generalized condition number}
$\kappa(A,B)$ is 
\[
\kappa(A,B)=\frac{\max\Lambda(A,B)}{\min\Lambda(A,B)}\;.
\]
We define the \emph{trace of $(A,B)$ }(denoted by $\trace(A,B)$)
as the sum of finite generalized eigenvalues of $(A,B$).
\end{defn}
(Generalized eigenvectors are defined also for matrices with different
null spaces~\cite{Stewart01}, but only the case of same null space
is relevant for this paper.) We will denote by $\Lambda(A)$ the set
of finite non-zero eigenvalues of $A$ (which is equal to $\Lambda(A,P_{A})$,
where $P_{A}$ is a projection onto the range of $A$). 

We are mainly interested in bounds on the smallest and largest generalized
eigenvalues (which we denote $\lambda_{\min}(\cdot,\cdot)$ and $\lambda_{\max}(\cdot,\cdot)$
respectively), since they tell us two important properties on the
pair $(A,B)$. First, for every unit norm vector $x$ we have 
\[
\lambda_{\min}(A,B)\cdot x^{T}Bx\leq x^{T}Ax\le\lambda_{\max}(A,B)\cdot x^{T}Bx\,.
\]
Second, when $B$ is used as a preconditioner for $A$, a vector $x$
satisfying $\Vert x-A^{+}b\Vert_{A}\leq\epsilon\Vert A^{+}b\Vert_{A}$
is found in at most $O(\sqrt{\kappa(A,B)}\cdot\log(1/\epsilon$))
iterations where $\Vert x\Vert_{A}^{2}=x^{T}Ax$ and $\kappa(A,B)=\lambda_{\max}(A,B)/\lambda_{\min}(A,B)$. 

In many cases it is easier to reason about non-generalized eigenvalues.
The following result from~\cite{ACST09} relates generalized eigenvalues
with regular eigenvalues of a different matrix.
\begin{lem}
\label{lemma:gev-to-singular}Let $A=UU^{T}$ and $B=VV^{T}$, where
$U$ and $V$ are real valued with the same number of rows. Assume
that $A$ and $B$ are symmetric, positive semidefinite and $\Null(A)=\Null(B)$.
We have
\[
\Lambda\left(A,B\right)=\Sigma^{2}\left(V^{+}U\right)
\]
and
\[
\Lambda\left(A,B\right)=\Sigma^{-2}\left(U^{+}V\right)\;.
\]
In these expressions, $\Sigma(\cdot)$ is the set of nonzero singular
values of the matrix within the parentheses, $\Sigma^{\ell}$ denotes
the same singular values to the $\ell$th power, and $V^{+}$ denotes
the Moore-Penrose pseudoinverse of $V$.
\end{lem}

\subsection{Effective resistance sampling}

Recent progress on fast SDD solvers~\cite{KoutisMillerPang10,BlellochEtAl11,KoutisMillerPang11}
is based on effective resistance sampling, first suggested in~\cite{SpielmanSrivastava08}.
Solving SDD systems can be reduced to solving a \emph{Laplacian }system.
Given a weighted undirected graph $G=([n],E,w)$ its \emph{Laplacian
}$L_{G}$ is given by $L=D-A$ where $A$ is the weighted adjacency
matrix $A_{ij}=w_{ij}$ and $D$ is the diagonal matrix of weighted
degrees given by $D_{ii}=\sum_{j\neq i}w_{ij}$. The \emph{effective
resistance }$R_{e}$ of an edge $e=(u,v)$ is given by 
\[
R_{e}=(e_{u}-e_{v})^{T}L_{G}^{+}(e_{u}-e_{v})
\]
where $e_{u}$ and $e_{v}$ are identity vectors and $L^{+}$ is the
Moore-Penrose pseudoinverse of $L$. The quantity is named effective
resistance because $R_{e}$ is equal to the potential difference induced
between $u$ and $v$ when a unit of current is injected at $u$ and
extracted at $v$, when $G$ is viewed as an electrical network with
conductances given by $w$.

Spielman and Srivastava~\cite{SpielmanSrivastava08} showed that
sampling sufficiently enough edges, where the probability of sampling
an edge is proportional to $w_{e}R_{e}$ , yields an high-quality
sparsifier $H$ for $G$. This implies that $L_{H}$ is an high-quality
preconditioner for $L_{G}$. Koutis et al.~\cite{KoutisMillerPang10,BlellochEtAl11,KoutisMillerPang11}
show that even crude approximations to the accurate effective resistances
suffice, and they show how such an approximation can be computed efficiently.
The asymptotically fastest solver~\cite{KoutisMillerPang11} solves
an $n$-by-$n$ SDD linear system in time $O(m\log n\log(1/\epsilon))$
where $m$ is the number of non-zeros in the matrix and $\epsilon$
is the accuracy of the solution.

\section{\label{sec:fem}Algebraic-Combinatorial Formulation of Finite Element
Matrices}

A finite element discretization of a PDE usually leads to an algebraic
system of equations $Kx=b$. The matrix $K$ has certain properties
that stem from the PDE and the specifics of how it was discretized.
To make our results more general and easier to understand by a wide
audience, we use the algebraic-combinatorial formulation developed
in~\cite{ShklarskiToledo08} rather than a PDE-derived formulation.

The matrix $K\in\mathbb{R}^{n\times n}$ is called a \emph{stiffness
matrix}, and it is a sum of \emph{element matrices}, $K=\sum_{e=1}^{m}K_{e}$.
Each element matrix $K_{e}$ corresponds to a subset of the domain
called a \emph{finite element}. The elements are disjoint except perhaps
for their boundaries and their union is the domain. We assume that
each element matrix $K_{e}$ is symmetric, positive semidefinite,
and zero outside a small set of $n_{e}$ rows and columns. In most
cases $n_{e}$ is uniformly bounded by a small integer. We denote
the set of nonzero rows and columns of $K_{e}$ by $\mathcal{N}_{e}$.
We denote the restriction of a matrix $A$ to indices $I$ by $A(I)$,
and denote the $\tilde{K}_{e}=K_{e}(\mathcal{N}_{e})$. $\tilde{K}_{e}$
is the \emph{essential element matrix }of $e$. Typically, in finite
element discretizations both the stiffness matrix ($K$) and the essential
element matrices ($\tilde{K}_{e}$s) are singular. For simplicity,
we assume that the rank and dimension of null space of all the elements
is the same and equal to $r$ and $d$ respectively. The null space
of $K$ is denoted by $\bm{\mathrm{N}}$ and we assume that its dimension
is $d$ as well.

Our proof technique relies on the fact that $K$ can be written as
$K=F^{T}F$ where 
\begin{equation}
F=\left(\begin{array}{c}
F_{1}\\
\vdots\\
F_{m}
\end{array}\right)\in\mathbb{R}^{mr\times n}\,.\label{eq:factored-form}
\end{equation}
In~\eqref{eq:factored-form} $F_{e}$ is the factored form $K_{e}$,
that is $K_{e}=F_{e}^{T}F_{e}$, so indeed $K=F^{T}F$. Many finite-element
discretization techniques actually generate the element matrices in
a factored form. Even if the elements are not generated in a factored
form, a factored form can be easily computed. One way to do so is
using the eigendecomposition $\tilde{K}_{e}=V_{e}\Sigma_{e}V_{e}^{T}$
. Define $\tilde{F}_{e}=\Sigma_{e}^{1/2}\bar{V}_{e}^{T}$ where $\bar{V}_{e}$
is obtained by taking the $r$ columns of $V_{e}$ associated with
non-zero eigenvalues, and let $F_{e}$ be obtained by expanding the
number of columns of $\tilde{F}_{e}$ to $n$ by adding zero columns
for columns not in $\mathcal{N}_{e}$. It is easy to verify that $K_{e}=F_{e}^{T}F_{e}$
and that $F_{e}$ is $r\times n$.

Typically, the factor has \emph{minimal rank deficiency }and the element
matrices are \emph{compatible }with $\bm{\mathrm{N}}$ and \emph{rigid}
with respect to it~\cite{ShklarskiToledo08}. We now explain what
these terms mean, as our theorems assumes that the finite element
discretization has them. We first discuss minimal rank deficiency.
\begin{defn}
A matrix $F\in\mathbb{R}^{m\times n}$ has\emph{ minimal rank deficiency}
if every set of $n-\dim(\Null(F))$ columns of $F$ is independent.

\end{defn}
Note that if the rank deficiency of $F$ is minimal then every leading
$l\times l$ minor of $K$ is non-singular, as long as $l\leq n-d$.
The null space $\bm{\mathrm{N}}$ of $K$ typically (that is, for
real-life finite element matrices) implies minimal rank deficiency,
but that has to be proven for each particular case. A simple technique
is based on the following lemma. 
\begin{lem}
\label{lem:minimal-technique}Suppose that $K=F^{T}F\in\mathbb{R}^{n\times n}$
has null space $\range(N)$ where $N\in\mathbb{R}^{n\times d}$. If
no $d\times d$ submatrix of $N$ is singular then $F$ has minimal
rank deficiency. \end{lem}
\begin{proof}
First notice that $\Null(F)=\Null(K)$ since $\Null(F^{T})=\range(F)^{\perp}$.
Suppose there is a set of $n-d$ columns of $F$ which are not independent.
Let $\bar{F}$ be a reordering of the columns of $F$ such that those
$n-d$ columns are first. There is a vector $x\in\mathbb{R}^{n-d}$
such that 
\[
\bar{F}\left(\begin{array}{c}
x\\
0_{d\times1}
\end{array}\right)=0\,.
\]
Let $\bar{N}$ be a reordering of the rows of $N$ consistently with
the reordering of the columns of $F$ in $\bar{F}$. The vector $\left(\begin{array}{cc}
x^{T} & 0\end{array}\right)^{T}$ is in the null space of $\bar{F}$ so there must exist a vector $y\neq0$
such that $\bar{N}y=\left(\begin{array}{cc}
x^{T} & 0\end{array}\right)^{T}$. This implies that the bottom $d$ rows of $\bar{N}$ form a singular
matrix. These rows are also rows of $N$, which implies that $N$
has a $d\times d$ singular submatrix, which contradicts our assumption.
\end{proof}
As an example, we show how Lemma~\ref{lem:minimal-technique} implies
minimal rank deficiency of the factor of a finite element matrix representing
a collection of elastic struts in two dimensions. In the next section
we show that Laplacians of connected graphs have minimal rank deficiency.
In~\cite{ShklarskiToledo08} it is shown that given a collection
$P=\left\{ p_{i}\right\} _{i=1}^{n}$ of points in the plane, the
null space of the rigid finite element matrix representing a collection
of elastic struts between the points is spanned by the range of 
\[
N=\left(\begin{array}{ccc}
1 & 0 & -y_{1}\\
0 & 1 & x_{1}\\
1 & 0 & -y_{2}\\
0 & 1 & x_{2}\\
\vdots & \vdots & \vdots\\
1 & 0 & -y_{n}\\
0 & 1 & x_{n}
\end{array}\right)\,.
\]
The matrix $N$ does not have singular $3$-by-$3$ submatrix unless
the points have some special properties (like three points with the
same $x$ coordinate), which they typically do not have. Even if such
a property is present, a slight rotation of the point set, an operation
that does not fundamentally change the physical problem, will remove
it.

We now turn to null space compatibility.
\begin{defn}
Let $A$ be an $m$-by-$n$ matrix, let $\mathcal{Z}_{A}$ be the
set of its zero columns. We define the \emph{essential null space
}of $A$ ($\enull(A)$) by 
\[
\enull(A)=\left\{ x\,:\, Ax=0\,\text{\,\,\ and}\,\,\, x_{i}=0\,\,\,\text{for}\,\,\, i\in\mathcal{Z}_{A}\right\} \,.
\]

\end{defn}

\begin{defn}
Let $\bm{\mathrm{N}}\subseteq\mathbb{R}^{n}$ be a linear space. A
matrix $A$ is called $\bm{\mathrm{N}}-$\emph{compatible }(or compatible
with $\bm{\mathrm{N}}$) if every vector in $x\in\enull(A)$ has a
unique vector $y\in\bm{\mathrm{N}}$ such that $x_{i}=y_{i}$ for
all $i\in\mathcal{N}_{A}$, and if the restriction of every vector
in $\bm{\mathrm{N}}$ to $\mathcal{N}_{A}$ (setting indices outside
$\mathcal{N}_{A}$ to zero) is always in $\enull(A)$.
\end{defn}
A particular discretization of a PDE yields element matrices ($K_{e}$s)
that are compatible with some well-known null space $\bm{\mathrm{N}}$,
which depends on the PDE; a translation in electrostatics, translations
and rotations in elasticity, and so on. Furthermore, it is usually
desirable that the stiffness matrix $K$ be \emph{rigid }with respect
to $\bm{\mathrm{N}}$, which is equivalent to saying that the null
space of $K$ is exactly $\bm{\mathrm{N}}$. For example, for matrix
of a resistive network elements are compatible with the span of the
all-ones vector. The null space of the the finite element matrix is
exactly the span of the all-ones (i.e., the matrix is rigid) if and
only if the graph is connected. Lack of rigidity often implies that
the PDE has not been discretized correctly, and it does not make sense
to solve the linear equations. This is an important scenario to detect
(see~\cite{ShklarskiToledo09}), but it is not the subject of this
paper. 

From now on we assume that the finite-element matrix $K$ has the
following \emph{well-formed} traits. 
\begin{lem}
The finite element matrix $K=F^{T}F$ is \emph{well-formed }if:
\begin{enumerate}
\item All elements are $\bm{\mathrm{N}}$-compatible.
\item $F$ has minimal rank deficiency.
\end{enumerate}
\end{lem}

\section{\label{sec:estiff}Effective Stiffness of an Element}

We now define the effective stiffness of an element. The stiffness
matrix of an element describes the physical properties (elasticity,
electrical conductivity, thermal conductivity, etc) of a piece of
material called an element by showing how that piece of material responds
to a load (current, mechanical force, etc) placed on the element.
The \emph{effective stiffness matrix} shows how the entire structure
responds to a load that is placed on one element. Intuitively, if
the stiffness matrix and the effective stiffness matrix of an element
are similar, the element is important; removing it from the structure
may significantly change the behavior of the overall structure. On
the other hand, if the effective stiffness element has a much larger
norm than the element matrix, then the element does not contribute
much to the strength (or conductivity) of the overall structure, so
it can be removed without changing much the overall behavior.

Algebraically, the effective stiffness matrix of $e$ is obtained
by eliminating (via Gauss elimination) from $K$ all columns not associated
with $e$.
\begin{defn}
\label{def:effective-stiffness}Assume that $K$ is well-formed. Let
$\bar{K}$ be obtained from $K$ by an arbitrary symmetric reordering
of the row and columns of $K$ such that the last $n_{e}$ rows and
columns of $\bar{K}$ are $\mathcal{N}_{e}$ and they are ordered
in ascending order (i.e., the ordering in $\bar{K}$ of the columns
in $\mathcal{N}_{e}$ is consistent with their order in $K$). Suppose
that $\bar{K}$ is partitioned 
\[
\bar{K}=\left(\begin{array}{cc}
\bar{K}_{11} & \bar{K}_{12}\\
\bar{K}_{12}^{T} & \bar{K}_{22}
\end{array}\right)
\]
 where $\bar{K}_{11}\in\mathbb{R}^{(n-n_{e})\times(n-n_{e})}$, $\bar{K}_{12}\in\mathbb{R}^{(n-n_{e})\times n_{e}}$
and $\bar{K}_{22}\in\mathbb{R}^{n_{e}\times n_{e}}$. The \emph{effective
stiffness} $S_{e}$ of element $e$ is 
\[
S_{e}=\bar{K}_{22}-\bar{K}_{12}^{T}\bar{K}_{11}^{-1}\bar{K}_{12}\,.
\]

\end{defn}
Note that the minimal rank deficiency of $K$ implies that $\bar{K}_{11}$
is non-singular, and that any ordering that respects the conditions
of the definition gives the same $S_{e}$, so the effective stiffness
is well defined. 

The following Lemma will be useful later on.
\begin{lem}
\label{lem:null}Assume that $K$ is well-formed. We have $\Null(S_{e})=\Null(\tilde{K}_{e})$
for every element $e$.\end{lem}
\begin{proof}
This lemma follows immediately from Lemma 3.7 and Lemma 5.5 from~\cite{ShklarskiToledo08}.
\end{proof}
Before proceeding to discuss effective stiffness sampling, and stating
our main result, we first show that indeed effective stiffness generalizes
effective resistance by showing that effective resistance is a particular
case of effective stiffness. 

The Laplacian of a weighted graph $G=([n],E,w)$ is, in fact, a finite
element matrix per our definition in section~\ref{sec:fem}. Given
an edge $e=(u,v)$ define $K_{e}=w_{e}(e_{u}-e_{v})(e_{u}-e_{v})^{T}$.
It is easy to verify that $L=\sum_{e\in E}K_{e}$. $L$ can also be
written in factor form $L=F^{T}F$ where $F\in\mathbb{R}^{\vert E\vert\times|V|}$.
Each edge $e=(u,v)$ correspond to row in $F$ given by $F_{e}=\sqrt{w_{e}}(e_{u}-e_{v})^{T}$.
It is well-known that if the if the graph is connected then the null
space of $L$ is exactly all-ones vector. Together with Lemma~\ref{lem:minimal-technique}
this implies that $F$ has minimal rank deficiency. It is also easy
to verify if $G$ that all elements are compatible with the all-ones
vector, so if $G$ is connected then $L$ is well-formed.

Simple calculation shows that $S_{e}1_{2}=0$ and $(e_{1}-e_{2})^{T}S_{e}(e_{1}-e_{2})=R_{e}^{-1}.$
This implies that $S_{e}=R_{e}^{-1}(e_{1}-e_{2})(e_{1}-e_{2})^{T}$
(here, $e_{1}=\left(\begin{array}{cc}
1 & 0\end{array}\right)^{T}$ and $e_{2}=\left(\begin{array}{cc}
0 & 1\end{array}\right)^{T}$).  Graph sparsification by effective resistance~\cite{SpielmanSrivastava08}
and near-linear time linear solvers~\cite{KoutisMillerPang10,BlellochEtAl11,KoutisMillerPang11}
relay on sampling edges with probability relative to $w_{e}R_{e}$.
It is easy to verify that $w_{e}R_{e}=\lambda_{\max}(\tilde{K}_{e},S_{e})$.
As we soon explain, we call the quantity $\lambda_{\max}(\tilde{K}_{e},S_{e})$
the \emph{leverage} of element $e$. Our main result shows that sampling
probabilities should be relative to the leverages for general finite
element matrices, and not only for Laplacians.

\section{Effective Stiffness Sampling}

This section defines the leverage of an of element and shows that
non-uniform sampling based on sampling probabilities that are relative
to the element leverages is a good choice. 
\begin{defn}
Assume that $K$ is well-formed. The \emph{leverage} of $e$ is 
\[
\tau_{e}=\lambda_{\max}(\tilde{K}_{e},S_{e})\,.
\]
(Recall that Lemma~\eqref{lem:null} guarantees that $\Null(S_{e})=\Null(\tilde{K}_{e})$).
The\emph{ total leverage} of $K$ is 
\[
\tau_{K}=\sum_{e=1}^{m}\tau_{e}\,.
\]
\end{defn}
\begin{note}
The term leverage arises from the connection between effective resistance
and statistical leverage that was noted by Drineas and Mahoney in~\cite{DrineasMahoney10}.
\end{note}
The main theorem shows how to use the leverages to sample finite element
matrices.
\begin{thm}
\label{thm:main}Let $K=F^{T}F=\sum_{e=1}^{m}K_{e}$ be an $n$-by-$n$
well-formed finite element matrix. Let 
\[
p_{e}=\frac{\tau_{e}}{\tau_{K}}
\]
and let $T_{1},\dots,T_{M}$ be a i.i.d random matrices defined by
\[
T_{i}=p_{J_{i}}^{-1}K_{J_{i}}
\]
where $J_{1},\dots,J_{M}$ are random integers between $1$ and $m$
which takes value $e$ with probability $p_{e}$. In other words,
$T_{i}$ is a scaled version of one of the $K_{e}$s, selected at
random, with a scaling that is proportional to the inverse of $p_{e}$.
Let $\kappa_{\max}>1$ and $\delta\in(0,1).$ If $M\geq C(\kappa_{\max})\tau_{K}\ln(2(n-d)/\delta)$
($C(\kappa_{\max})$ is given by~\eqref{eq:Ceq}) then 
\[
\Pr\left(\Null\left(\frac{1}{M}\sum_{i=1}^{M}T_{i}\right)\neq\bm{\mathrm{N}}\,\,\mathrm{or}\,\,\kappa\left(K,\frac{1}{M}\sum_{i=1}^{M}T_{i}\right)>\kappa_{\max}\right)\leq\delta\,.
\]

\end{thm}
Before proving Theorem~\ref{thm:main} we need to state and prove
a few auxiliary lemmas. In the following two lemmas, $K=F^{T}F=\sum_{e=1}^{m}K_{e}$
is a $n$-by-$n$ well-formed finite element matrix.\label{lem:probs}Let
$U\in\mathbb{R}^{mr\times n}$ be any matrix whose columns form an
orthonormal basis of $\range(F)$. Let $U_{e}\in\mathbb{R}^{r\times n}$
be the rows of $U$ corresponding to element $e$. The set of non-zero
eigenvalues (including multiplicity) of $U_{e}U_{e}^{T}$ and the
set of finite generalized eigenvalues of $(\tilde{K}_{e},S_{e})$
are the same. In particular, 
\[
\lambda_{\max}(U_{e}U_{e}^{T})=\lambda_{\max}(\tilde{K}_{e},S_{e})=\tau_{e}
\]
and 
\[
\trace(U_{e}U_{e}^{T})=\trace(\tilde{K}_{e},S_{e})\,.
\]

\begin{proof}
We first show that we can prove the lemma by showing that it holds
for a particular $U$. An arbitrary orthonormal basis $V$ is related
to $U$ by $V=UZ$, where $Z$ is an $n$-by-$n$ unitary matrix.
In particular, $V_{e}=U_{e}Z$ ($V_{e}$ are the rows of $V$ corresponding
to element $e$) so $V_{e}V_{e}^{T}=U_{e}ZZ^{T}U_{e}^{T}=U_{e}U_{e}^{T}$.
We obtain $U$ from the $QR$ factorization of $\bar{F}=\bar{U}R$
and set $U$ to be the first $n-d$ columns of $\bar{U}$, where $\bar{F}$
is obtained from $F$ by reordering the columns in $\mathcal{N}_{e}$
to the end (consistently with their ordering in $F$).

The last $n_{e}$ columns of $\bar{F}$ are $\mathcal{N}_{e}$, and
$F_{e}$ is non-zero outside the indices of $\mathcal{N}_{e}$.This
implies that 
\[
\bar{F}_{e}=\left[\begin{array}{cc}
0_{r\times(n-n_{e})} & \tilde{F}_{e}\end{array}\right]
\]
\[
U_{e}=\left[\begin{array}{cc}
0_{r\times(n-n_{e})} & \tilde{U}_{e}\end{array}\right]
\]
where $\tilde{U}_{e},\,\tilde{F}_{e}\in\mathbb{R}^{r\times n_{e}}$.
Let us write 
\[
R=\left(\begin{array}{cc}
R_{11} & R_{12}\\
0 & R_{22}
\end{array}\right)
\]
where $R_{11}\in\mathbb{R}^{(n-n_{e})\times(n-n_{e})}$, $R_{12}\in\mathbb{R}^{(n-n_{e})\times n_{e}}$
and $R_{22}\in\mathbb{R}^{n_{e}\times n_{e}}$. Let us write $\bar{K}=\bar{F}^{T}\bar{F}$
and 
\[
\bar{K}=\left(\begin{array}{cc}
\bar{K}_{11} & \bar{K}_{12}\\
\bar{K}_{12}^{T} & \bar{K}_{22}
\end{array}\right)
\]
where $\bar{K}_{11}\in\mathbb{R}^{(n-n_{e})\times(n-n_{e})}$, $\bar{K}_{12}\in\mathbb{R}^{(n-n_{e})\times n_{e}}$
and $\bar{K}_{22}\in\mathbb{R}^{n_{e}\times n_{e}}$. Since $R$ is
the $R$-factor of $\bar{F}$ and $\bar{K}=\bar{F}^{T}\bar{F}$ it
is also the Cholesky factor of $\bar{K}$. It also implies that $R_{22}^{T}R_{22}$
is equal to the Schur complement  
\[
R_{22}^{T}R_{22}=\bar{K}_{22}-\bar{K}_{12}^{T}\bar{K}_{11}^{-1}\bar{K}_{12}=S_{e}\,.
\]

The minimal rank deficiency of $F$ implies that that the bottom $d$
rows of $R$ and $R_{22}$ are zero. Let $\bar{R}_{22}\in\mathbb{R}^{(n_{e}-d)\times n_{e}}$
be the first $n_{e}-d$ rows of $R_{22}$. It is still the case that
$\bar{R}_{22}^{T}\bar{R}_{22}=S_{e}$. We have $\bar{F}=\bar{U}R$,
so $\bar{F}_{e}=\bar{U}_{e}R_{22}=U_{e}\bar{R}_{22}$ which implies
that $\tilde{F}_{e}=\tilde{U}_{e}\bar{R}_{22}$. Applying Lemma~\ref{lemma:gev-to-singular}
we find that 
\begin{eqnarray*}
\Lambda(\tilde{K}_{e},S_{e}) & = & \Lambda(\tilde{F}_{e}^{T}\tilde{F}_{e},R_{22}^{T}R_{22})\\
 & = & \Sigma^{2}\left((\bar{R}_{22}^{T})^{+}\tilde{F}_{e}^{T}\right)\\
 & = & \Sigma^{2}\left((\bar{R}_{22}^{T})^{+}\bar{R}_{22}^{T}\tilde{U}_{e}^{T}\right)
\end{eqnarray*}
The minimal rank deficiency of $\bar{F}$ implies that $R_{22}$ is
full rank, so $R_{22}^{T}$ is a full rank matrix with more rows than
columns (or equal), so $(\bar{R}_{22}^{T})^{+}\bar{R}_{22}^{T}=I_{n_{e}}$.
This implies that $(\bar{R}_{22}^{T})^{+}\bar{R}_{22}^{T}\tilde{U}_{e}^{T}=\tilde{U}_{e}^{T}$
so
\[
\Lambda(\tilde{K}_{e},S_{e})=\Sigma^{2}(\tilde{U}_{e}^{T})\,.
\]
$\Sigma^{2}(\tilde{U}_{e}^{T})$ is exactly the set of non-zero eigenvalues
of $\tilde{U}_{e}\tilde{U}_{e}^{T}$. Therefore, the non-zero eigenvalues
of $U_{e}U_{e}^{T}$ are exactly the finite generalized eigenvalues
of $(\tilde{K}_{e},S_{e})$, so 
\[
\lambda_{\max}(U_{e}U_{e}^{T})=\lambda_{\max}(\tilde{K}_{e},S_{e})=\tau_{e}
\]
 and 
\[
\trace(U_{e}U_{e}^{T})=\trace(\tilde{K}_{e},S_{e})\,.
\]
\end{proof}
\begin{lem}
\label{lem:leverage-bound} We have $(n-d)/r\leq\tau_{K}\leq n-d$.\end{lem}
\begin{proof}
\begin{eqnarray*}
\tau_{K}=\sum_{e=1}^{m}\tau_{e}=\sum_{e=1}^{m}\lambda_{\max}(\tilde{K}_{e},S_{e}) & \leq & \sum_{e=1}^{m}\trace(\tilde{K}_{e},S_{e})\\
 & = & \sum_{e=1}^{m}\trace(U_{e}U_{e}^{T})\\
 & = & \sum_{e=1}^{m}\trace(U_{e}^{T}U_{e})\\
 & = & \trace(\sum_{i=1}^{m}U_{e}^{T}U_{e})\\
 & = & \trace(U^{T}U)=n-d\,.
\end{eqnarray*}
For each element the pencil $(\tilde{K}_{e},S_{e})$ has exactly $r$
determined eigenvalues, so $\lambda_{\max}(\tilde{K}_{e},S_{e})\geq\trace(\tilde{K}_{e},S_{e})/r$.
The lower bound follows.
\end{proof}
We can now prove Theorem~\ref{thm:main}.
\begin{proof}
\emph{(of Theorem~\ref{thm:main}) }We express the matrix $\frac{1}{M}\sum_{i=1}^{M}T_{i}$
as a normal form\emph{ }
\begin{equation}
\frac{1}{M}\sum_{i=1}^{M}T_{i}=(\mathcal{S}F)^{T}(\mathcal{S}F)\label{eq:factor-form}
\end{equation}
where $\mathcal{S}\in\mathbb{R}^{Mr\times mr}$ is a random sampling
matrix and $F$ is the factor of the stiffness matrix $K=F^{T}F$.
If we take ${\cal S}$ to be a block matrix with $r\times r$ blocks,
its blocks defined by

\[
\mathcal{S}_{ie}=\begin{cases}
\sqrt{\frac{1}{M}}p_{e}^{-1/2}I_{r\times r} & \text{if }T_{i}=p_{e}^{-1}K_{e}\\
0_{r\times r} & \text{otherwise,}
\end{cases}
\]
then it is easy to verify that~equation~\eqref{eq:factor-form}
is satisfied. Let $F=\bar{U}\bar{R}$ be a reduced $QR$ factorization
of $F$. The minimal rank deficiency of $F$ implies that that the
bottom $d$ rows of $R$ are zero. Let $R\in\mathbb{R}^{(n-d)\times n}$
be the first $n-d$ rows of $\bar{R}$, and $U\in\mathbb{R}^{mr\times(n-d)}$
be the first $n-d$ columns of $\bar{U}$. It is easy to verify that
$F=UR$ and $F^{T}F=R^{T}R$. $R^{T}$ is full rank, so $(R^{T})^{+}R^{T}=I_{n}$.
Assume for now that $\Null(\mathcal{S}F)=\Null(F)$. Applying lemma~\ref{lemma:gev-to-singular}
we have 
\begin{eqnarray*}
\kappa(K,\frac{1}{M}\sum_{i=1}^{M}T_{i}) & = & \kappa(F^{T}F,(\mathcal{S}F)^{T}(\mathcal{S}F))\\
 & = & \kappa(R^{T}R,(\mathcal{S}F)^{T}(\mathcal{S}F))\\
 & = & \kappa^{2}\left((R^{T})^{+}F^{T}\mathcal{S}^{T}\right)\\
 & = & \kappa^{2}\left((R^{T})^{+}R^{T}U^{T}\mathcal{S}^{T}\right)\\
 & = & \kappa^{2}(U^{T}\mathcal{S}^{T})\\
 & = & \kappa^{2}((\mathcal{S}U)^{T})\\
 & = & \kappa^{2}(\mathcal{S}U)\\
 & = & \kappa((\mathcal{S}U)^{T}(\mathcal{S}U))\,.
\end{eqnarray*}

Define the i.i.d random matrices $Y_{1},\dots,Y_{M}$ by 
\[
Y_{i}=p_{J_{i}}^{-1}U_{J_{i}}^{T}U_{J_{i}}
\]
where $U_{e}$ is the rows corresponding to element $e$ in $U$.
It is easy to verify that 
\[
(\mathcal{S}U)^{T}(\mathcal{S}U)=\frac{1}{M}\sum_{i=1}^{M}Y_{i}\,.
\]

If $\Null(\mathcal{S}F)=\Null(F)$ then $\Null(\frac{1}{M}\sum_{i=1}^{M}T_{i})=\bm{\mathrm{N}}$.
$U$ is full rank so $\Null(F)=\Null(UR)=\Null(R)$. On the other
hand $\mathcal{S}F=\mathcal{S}UR$, so $\mathcal{S}U$ is full rank
if and only if $\Null(\mathcal{S}F)=\Null(R)=\Null(F)$. $\mathcal{S}U$
is rank deficient only if $\frac{1}{M}\sum_{i=1}^{M}Y_{i}$ is singular.
Furthermore, if $\frac{1}{M}\sum_{i=1}^{M}Y_{i}$ is not singular,
then $\Null(\mathcal{S}F)=\Null(F)$ as required earlier.

Combining previous arguments, we find that 
\[
\Pr\left(\Null\left(\frac{1}{M}\sum_{i=1}^{M}T_{i}\right)\neq\bm{\mathrm{N}}\,\,\mathrm{or}\,\,\kappa\left(K,\frac{1}{M}\sum_{i=1}^{M}T_{i}\right)>\kappa_{\max}\right)\leq\Pr\left(\frac{1}{M}\sum_{i=1}^{M}Y_{i}\,\mbox{\,\textrm{is singular}}\,\,\mathrm{or}\,\,\kappa\left(\frac{1}{M}\sum_{i=1}^{M}Y_{i}\right)>\kappa_{\max}\right)\,.
\]
The expectation of the $Y_{i}$'s is the identity matrix, 
\begin{eqnarray*}
\E(Y_{i}) & = & \sum_{j=1}^{M}\Pr(T_{i}=p_{j}^{-1}K_{j})p_{j}^{-1}U_{j}^{T}U_{j}\\
 & = & \sum_{j=1}^{M}p_{j}p_{j}^{-1}U_{j}^{T}U_{j}\\
 & = & \sum_{j=1}^{M}U_{j}^{T}U_{j}\\
 & = & U^{T}U=I_{n-d\times n-d}
\end{eqnarray*}
and their $2$-norm is bounded by 
\begin{eqnarray*}
\left\Vert Y_{i}\right\Vert _{2} & \leq & \max_{j}p_{j}^{-1}\left\Vert U_{j}^{T}U_{j}\right\Vert _{2}\\
 & = & \max_{j}p_{j}^{-1}\lambda_{\max}(U_{j}U_{j}^{T})\\
 & = & \max_{j}p_{j}^{-1}\lambda_{\max}(\tilde{K}_{j},S_{j})\\
 & = & \max_{j}p_{j}^{-1}\tau_{j}\\
 & = & \max_{j}\left(\left(\frac{\tau_{j}}{\tau_{K}}\right)^{-1}\tau_{j}\right)\\
 & = & \tau_{K}\leq n-d.
\end{eqnarray*}

We now apply Corollary~\eqref{cor:kappa-bound} on $Y_{1},\dots,Y_{M}$
to find that 
\[
\Pr\left(\frac{1}{M}\sum_{i=1}^{M}Y_{i}\,\,\mbox{\textrm{is singular}}\,\,\mathrm{or}\,\,\kappa\left(\frac{1}{M}\sum_{i=1}^{M}Y_{i}\right)>\kappa_{\max}\right)\leq\delta
\]

\end{proof}

\paragraph{Comparison to Spielman and Srivastava's bound for effective resistance
sampling~\cite{SpielmanSrivastava08}.}

Effective resistance sampling is a case of effective stiffness sampling.
If we examine the sampling procedure analyzed in \cite[Theorem~1]{SpielmanSrivastava08}
we see that for Laplacians it is identical to the the one analyzed
in Theorem~\ref{thm:main}. We now compare the analyses. 

There are two differences in the way the bounds are formulated: 
\begin{enumerate}
\item Spielman and Srivastava are mainly interested in spectral partitioning,
so they compare the sparsified quadratic form $x^{T}L_{H}x$ to the
original quadratic form $x^{T}L_{G}x$. We are mainly interested in
using the sparified matrix as a preconditioner, so we bound the maximum
condition number $\kappa_{\max}$. However, it is easy to modify our
analysis to give bounds in terms of quadratic forms. On the other
hand, Spielman and Srivastava's bound immediately leads to a $(1+\epsilon)/(1-\epsilon)$
bound on the condition number. Using $\epsilon=(\kappa_{\max}-1)/(\kappa_{\min}+1)$
we can convert Spielman and Srivastava's bound to a bound in terms
of $\kappa_{\max}$.
\item Spielman and Srivastava's bound applies only for one failure probability:
$1/2$ (however, the analysis might be modified to allow other failure
probabilities). 
\end{enumerate}
Setting $\delta=1/2$, our bound for Laplacians ($d=1,\,\tau_{K}=n-1$)
for this failure probability is $M\geq C(\kappa_{\max})(n-1)\ln(4(n-1))$\@.
Writing $\epsilon=(\kappa_{\max}-1)/(\kappa_{\max}+1)$ we find that
Spielman and Srivastava bound is $\tilde{M}\geq\tilde{C}(\kappa_{\max})n\ln(n)$
where 
\[
\tilde{C}(\kappa_{\max})=\frac{9(\kappa_{\max}+1)^{2}R}{(\kappa_{\max}-1)^{2}}\,.
\]
$R$ is some unspecified constant. The unspecified constant $R$ makes
a comparison hard (and might also cause problems trying to apply the
theorem). However, if assume $R=1$, then $\tilde{C}(3)=36$ and $\lim\tilde{C}(\kappa_{\max})=9$
(where $\kappa_{\max}=3$ is taken as an example). The constants in
Theorem~\eqref{thm:main} are $C(3)\approx9.2423$ and $\lim C(\kappa_{\max})=1/(2\ln2-1)\thickapprox2.5887$.
So, it seems that the constants in our bound are much better, although
this is probably partly due to the fact that we are using newer and
tighter matrix Chernoff bounds (Theorem~\ref{thm:Chernoff}). Asymptotically,
both bounds are equivalent.

\section{Sampling Using Inexact Leverages or Upper Bounds}

Theorem~\ref{thm:main} shows that the sampling probabilities that
are proportional to $\tau_{e}$ are effective for randomly selecting
a good subset of elements to serve as a preconditioner. In practice
it may be possible to obtain only estimates for the true maximum eigenvalues.
The following two generalizations of Theorem~\ref{thm:main} show
that even crude approximations or upper bounds of the leverages suffice,
provided that the number of samples is enlarged accordingly.
\begin{thm}
\label{thm:approx1}For every element $e$ let $\tilde{\tau}_{e}$
be $(1+\delta)$-approximations to $\tau_{e}$, that is 
\[
\left|\tilde{\tau}_{e}-\tau_{e}\right|\leq\delta\cdot\tau_{e}\,.
\]
We make the same assumptions and use the same notation as in Theorem~\ref{thm:main}
except that the probabilities $p_{e}$ are now given by 
\[
p_{e}=\frac{\tilde{\tau}_{e}}{\sum_{i=1}^{m}\tilde{\tau}_{i}}\,.
\]
 If $M\geq C(\kappa_{\max})\tau_{K}\beta\ln(2(n-d)/\delta)$ ($C(\kappa_{\max})$
is given by~\eqref{eq:Ceq}) , where $\beta=\frac{1+\delta}{1-\delta}$,
then 
\[
\Pr\left(\Null\left(\frac{1}{M}\sum_{i=1}^{M}T_{i}\right)\neq\bm{\mathrm{N}}\,\,\mathrm{or}\,\,\kappa\left(K,\frac{1}{M}\sum_{i=1}^{M}T_{i}\right)>\kappa_{\max}\right)\leq\delta\,.
\]
\end{thm}
\begin{proof}
The proof is identical to the proof of Theorem~\ref{thm:main} except
that the bound on $\Vert Y_{i}\Vert_{2}$ needs to be modified as
follows: 
\begin{eqnarray*}
\left\Vert Y_{i}\right\Vert _{2} & \leq & \max_{j}p_{j}^{-1}\left\Vert U_{j}^{T}U_{j}\right\Vert _{2}\\
 & = & \max_{j}p_{j}^{-1}\lambda_{\max}(U_{j}U_{j}^{T})\\
 & = & \max_{j}p_{j}^{-1}\lambda_{\max}(\tilde{K}_{j},S_{j})\\
 & = & \max_{j}\left(\left(\frac{\tilde{\tau}_{j}}{\sum_{i=1}^{m}\tilde{\tau}_{i}}\right)^{-1}\tau_{j}\right)\\
 & \leq & \max_{j}\left(\left(\frac{(1-\delta)\tau_{j}}{(1+\delta)\sum_{i=1}^{m}\tau_{i}}\right)^{-1}\tau_{j}\right)\\
 & = & \beta\sum_{e=1}^{m}\tau_{e}\\
 & \leq & \tau_{K}\beta\,.
\end{eqnarray*}
\end{proof}
\begin{thm}
\label{thm:approx2}For every element $e$ let $\tilde{\tau}_{e}$
be and upper bound on $\tau_{e}$, and let $\tilde{\tau}_{K}=\sum_{e=1}^{m}\tilde{\tau}_{e}$.
We make the same assumptions and use the same notation as in Theorem~\ref{thm:main}
except that the probabilities $p_{e}$ are now given by 
\[
p_{e}=\frac{\tilde{\tau}_{e}}{\tilde{\tau}_{K}}\,.
\]
 If $M\geq C(\kappa_{\max})\tilde{\tau}_{K}\ln(2(n-d)/\delta)$ ($C(\kappa_{\max})$
is given by~\eqref{eq:Ceq}) then 
\[
\Pr\left(\Null\left(\frac{1}{M}\sum_{i=1}^{M}T_{i}\right)\neq\bm{\mathrm{N}}\,\,\mathrm{or}\,\,\kappa\left(K,\frac{1}{M}\sum_{i=1}^{M}T_{i}\right)>\kappa_{\max}\right)\leq\delta\,.
\]
\end{thm}
\begin{proof}
The proof is identical to the proof of Theorem~\ref{thm:main} except
that the bound on $\Vert Y_{i}\Vert_{2}$ needs to be modified as
follows: 
\begin{eqnarray*}
\left\Vert Y_{i}\right\Vert _{2} & \leq & \max_{j}p_{j}^{-1}\left\Vert U_{j}^{T}U_{j}\right\Vert _{2}\\
 & = & \max_{j}p_{j}^{-1}\lambda_{\max}(U_{j}U_{j}^{T})\\
 & = & \max_{j}p_{j}^{-1}\lambda_{\max}(\tilde{K}_{j},S_{j})\\
 & = & \max_{j}\left(\left(\frac{\tilde{\tau}_{j}}{\tilde{\tau}_{K}}\right)^{-1}\tau_{j}\right)\\
 & \leq & \tilde{\tau}_{K}\cdot\max_{j}\frac{\tau_{j}}{\tilde{\tau}_{j}}\\
 & \leq & \tilde{\tau}_{K}\,.
\end{eqnarray*}

\end{proof}

\section{A Condition-number Formula for The Leverages}

In this section we show that the leverage $\tau_{e}$ can also be
defined in terms of the condition number of $(K,K-K_{e})$. This condition
number is the one related to preconditioning $K$ by removing only
element $e$. 
\begin{thm}
Let $K=F^{T}F=\sum_{e=1}^{m}K_{e}$ be an $n$-by-$n$ well-formed
finite element matrix. For every element $e$, if $\Null(K-K_{e})=\Null(K)$
then 
\[
\tau_{e}=\frac{\kappa(K,K-K_{e})-1}{\kappa(K,K-K_{e})}\,,
\]
 otherwise $\tau_{e}=1$. \end{thm}
\begin{proof}
We now argue that if $\rank(K-K_{e})<\rank(K)$ then $\tau_{e}=1$.
Let $\bar{K}$ be obtained from $K-K_{e}$ by an arbitrary symmetric
reordering of the row and columns of $K$ such that the last $n_{e}$
rows and columns of $\bar{K}$ are $\mathcal{N}_{e}$ and they are
ordered in ascending order (i.e., the ordering in $\bar{K}$ of the
columns in $\mathcal{N}_{e}$ is consistent with their order in $K$).
Suppose that $\bar{K}$ is partitioned 
\[
\bar{K}=\left(\begin{array}{cc}
\bar{K}_{11} & \bar{K}_{12}\\
\bar{K}_{12}^{T} & \bar{K}_{22}
\end{array}\right)
\]
 where $\bar{K}_{1}\in\mathbb{R}^{(n-n_{e})\times(n-n_{e})}$, $\bar{K}_{12}\in\mathbb{R}^{(n-n_{e})\times n_{e}}$
and $\bar{K}_{22}\in\mathbb{R}^{n_{e}\times n_{e}}$. is well-formed
so $\bar{K}_{11}$ is non-singular. This implies that $\rank(K-K_{e})=\rank(\bar{K})=n-n_{e}+\rank(\bar{K}_{22}-\bar{K}_{12}^{T}\bar{K}_{11}^{-1}\bar{K}_{12})$
since $\bar{K}_{22}-\bar{K}_{12}^{T}\bar{K}_{11}^{-1}\bar{K}_{12}$
is the Schur complement. It is easy to see that $\bar{K}_{22}-\bar{K}_{12}^{T}\bar{K}_{11}^{-1}\bar{K}_{12}=S_{e}-\tilde{K}_{e}$.
On the other hand, using similar observations we find that $\rank(K)=n-n_{e}+\rank(S_{e})$.
From $\rank(K-K_{e})<\rank(K)$ we find that $\rank(S_{e}-\tilde{K}_{e})<\rank(S_{e})$.
Therefore there exists a vector $x$ such that $S_{e}x\neq0$ but
$(S_{e}-\tilde{K}_{e})x=0$. That $x$ is an eigenvector of $(\tilde{K}_{e},S_{e})$
corresponding to the eigenvalue $1$ since we have $\tilde{K}_{e}x=S_{e}x$
but $S_{e}x\neq0$. All eigenvalues of $(\tilde{K}_{e},S_{e})$ are
bounded by $1$ so we found that $\lambda_{\max}(\tilde{K}_{e},S_{e})=1$.

We now analyze the spectrum of $(K,K-K_{e})$. Without loss of generality
assume $e=m$. Let $\mathcal{S}\in\mathbb{R}^{(m-1)r\times mr}$ be
defined as 
\[
\mathcal{S}=\left[\begin{array}{cc}
I_{(m-1)r\times mr} & 0_{(m-1)r\times r}\end{array}\right]\,.
\]
It is easy to verify that $K-K_{e}=(\mathcal{S}F)^{T}(\mathcal{S}F)$. 

Let $F=\bar{U}\bar{R}$ be a reduced $QR$ factorization of $F$.
The minimal rank deficiency of $F$ implies that that the bottom $d$
rows of $\bar{R}$ are zero. Let $R\in\mathbb{R}^{(n-d)\times n}$
be the first $n-d$ rows of $\bar{R}$, and $U\in\mathbb{R}^{mr\times(n-d)}$
be the first $n-d$ columns of $\bar{U}$. It is easy to verify that
$F=UR$ and $F^{T}F=R^{T}R$. The matrix $R^{T}$ is full rank, so
$(R^{T})^{+}R^{T}=I_{n}$. $U$ has orthonormal rows so $U^{T}U=I_{(n-d)\times(n-d)}$.
Applying lemma~\ref{lemma:gev-to-singular} we have 
\begin{eqnarray*}
\Lambda(K,K-K_{e}) & = & \Lambda(F^{T}F,(\mathcal{S}F)^{T}(\mathcal{S}F))\\
 & = & \Lambda(R^{T}R,(\mathcal{S}F)^{T}(\mathcal{S}F))\\
 & = & \Sigma^{2}\left((R^{T})^{+}F^{T}S^{T}\right)\\
 & = & \Sigma^{2}\left((R^{T})^{+}R^{T}U^{T}S^{T}\right)\\
 & = & \Sigma^{2}(U^{T}S^{T})\\
 & = & \Lambda(U^{T}\mathcal{S}^{T}\mathcal{S}U)\,.
\end{eqnarray*}

Let $T\in\mathbb{R}^{mr\times mr}$ be defined as 
\[
T=\left[\begin{array}{cc}
0_{(m-1)\times(m-1)r}\\
 & I_{r\times r}
\end{array}\right]\,.
\]
It is easy to verify that $\mathcal{S}^{T}\mathcal{S}=I_{mr\times mr}-T$~.
We now have 
\begin{eqnarray*}
\Lambda(U^{T}\mathcal{S}^{T}\mathcal{S}U) & = & \Lambda(U^{T}(I_{mr\times mr}-T)U)\\
 & = & \Lambda(U^{T}U-U^{T}TU)\\
 & = & \Lambda(I_{(n-d)\times(n-d)}-U^{T}TU)\,.
\end{eqnarray*}
Let $U_{e}$ be the bottom $r$ rows of $U$. It is easy to verify
that $U^{T}TU=U_{e}^{T}U_{e}$, so $\Lambda(U^{T}\mathcal{S}^{T}\mathcal{S}U)=\Lambda(I-U_{e}^{T}U_{e})$.
Let $(\lambda,x)$ be an eigenpair of $U_{e}^{T}U_{e}$, that is $U_{e}^{T}U_{e}x=\lambda x$.
We have 
\[
(I-U_{e}^{T}U_{e})x=x-U_{e}^{T}U_{e}x=x-\lambda x=(1-\lambda)x\,,
\]
so $(1-\lambda,x)$ is an eigenpair of $I-U_{e}^{T}U_{e}$. $U_{e}^{T}U_{e}$
is an order $n-d$ matrix of rank $r<n-d$ so it is singular. $U_{e}^{T}U_{e}$
is also positive semidefinite so all its eigenvalues are non-negative.
The last three facts imply that $\lambda_{\max}(I-U_{e}^{T}U_{e})=1$.
On the other hand, clearly $\lambda_{\min}(I-U_{e}^{T}U_{e})=1-\lambda_{\max}(U_{e}^{T}U_{e})$.
Combining these two together we find that 
\[
\kappa(K,K-K_{e})=\kappa(I-U_{e}^{T}U_{e})=\frac{1}{1-\lambda_{\max}(U_{e}^{T}U_{e})}\,.
\]

This implies that 
\[
\frac{\kappa(K,K-K_{e})-1}{\kappa(K,K-K_{e})}=\lambda_{\max}(U_{e}^{T}U_{e})\,.
\]
The non-zero eigenvalues of $U_{e}U_{e}^{T}$ are exactly the non-zero
eigenvalues of $U_{e}^{T}U_{e}$, so $\lambda_{\max}(U_{e}U_{e}^{T})=\lambda_{\max}(U_{e}^{T}U_{e})$.
$U$ is a matrix whose columns form an orthonormal basis of $\range(F)$,
so according to Lemma~\ref{lem:probs} we have $\lambda_{\max}(U_{e}U_{e}^{T})=\tau_{e}$,
which concludes the proof.
\end{proof}

\section{Rayleigh Monotonicity Law for Finite Element Matrices and Local Approximation
of Effective Stiffness}

For electrical circuits it is well known that when the resistances
of a circuit are increased, the effective resistance between any two
points can only increase. If the resistances are decreased, the effective
resistance can only decrease. This is the so-called ``Rayleigh Monotonicity
Law''. The following theorem shows that a similar statement can be
said about the effective stiffness.
\begin{thm}
[Rayleigh Monotonicity Law for Finite Element Matrices]\label{thm:rayleigh}
Let $K=F^{T}F=\sum_{e=1}^{m}K_{e}$ be an $n$-by-$n$ well-formed
finite element matrix. Assume that for every element $e$ we have
a factorization $K_{e}=B_{e}^{T}R_{e}^{-1}B_{e}$ such that $R_{e}\in\mathbb{R}^{r\times r}$
is symmetric positive definite and $B_{e}\in\mathbb{R}^{r\times n}$
has rank $r$. Let $\hat{K}_{e}=\sum_{e=1}^{m}\hat{K}_{e}$ be another
finite element matrix with the same set of non-zero rows and columns
for every element $e$, and assume every element has a factorization
$\hat{K}_{e}=B_{e}^{T}\hat{R}_{e}^{-1}B_{e}$ such that $\hat{R}_{e}\in\mathbb{R}^{r\times r}$
is symmetric positive definite, and $R_{e}\preceq\hat{R}_{e}$ for
every $e$. For an element $e$, let $S_{e}$ be the effective stiffness
of $e$ in $K$, and $\hat{S}_{e}$ be the effective stiffness of
$e$ in $\hat{K}$. Then $S_{e}^{+}\preceq\hat{S}_{e}^{+}$. \end{thm}
\begin{proof}
Denote
\[
B=\left[\begin{array}{c}
B_{1}\\
\vdots\\
B_{m}
\end{array}\right],\,\,\, R=\left(\begin{array}{cccc}
R_{1}\\
 & R_{2}\\
 &  & \ddots\\
 &  &  & R_{m}
\end{array}\right),\,\text{and}\,\hat{R}=\left(\begin{array}{cccc}
\hat{R}_{1}\\
 & \hat{R}_{2}\\
 &  & \ddots\\
 &  &  & \hat{R}_{m}
\end{array}\right)\,.
\]
Notice that $K=B^{T}R^{-1}B$, $\hat{K}=B^{T}\hat{R}^{-1}B$, and
that $R\preceq\hat{R}$.

Since $B_{e}$ has full row rank and both $R_{e}$ and $\hat{R}_{e}$
are non-singular, we have $\Null(K_{e})=\Null(\hat{K}_{e}$). This
implies that $K$ and $\hat{K}$ are compatible with the same null
space $\bm{\mathrm{N}}$ . This, in turn, implies that $\Null(S_{e})=\Null(\hat{S}_{e})$
(Lemma~\ref{lem:null}), so it is enough to prove that for every
$x\perp\Null(S_{e})$ we have $x^{T}S_{e}^{+}x\leq x^{T}\hat{S}_{e}^{+}x$
. 

Fix some element $e$. To avoid notation clutter we will assume, without
loss of generality, that $K$ and $\hat{K}$ are ordered as in Definition~\ref{def:effective-stiffness}.
Let $x\perp\Null(S_{e})$ and let $y=(\begin{array}{cc}
0_{1\times(n-n_{e})} & x\end{array})^{T}.$ It is easy to verify that $x^{T}S_{e}^{+}x=y^{T}K^{+}y$. Let $f=R^{-1}BK^{+}y$.
We now have 
\[
f^{T}Rf=y^{T}K^{+}B^{T}R^{-1}RR^{-1}K^{+}y=y^{T}K^{+}KK^{+}y=y^{T}K^{+}y
\]
where the last equality follows since $y\perp\Null(K)$ ($\bm{\mathrm{N}}$-compatibility).

Since $R\preceq\hat{R}$ we have $f^{T}Rf\leq f^{T}\hat{R}f$. Let
\begin{equation}
\hat{f}=\arg\min_{B^{T}g=y}g^{T}\hat{R}g\,.\label{eq:primal}
\end{equation}
Since $B^{T}f=KK^{+}y=y$ we have $f^{T}\hat{R}f\leq\hat{f}^{T}\hat{R}\hat{f}$.
We now have $\hat{f}^{T}\hat{R}\hat{f}=y^{T}\hat{K}^{+}y$ since the
minimization~\eqref{eq:primal} is dual to $\max_{v\in\mathbb{R}^{n}}2v^{T}y-v^{T}\hat{K}v$
whose maximum is attained at $v=\hat{K}^{+}y$. Finally, it is easy
to verify that $y^{T}\hat{K}^{+}y=x^{T}\hat{S}_{e}^{+}x$. Combining
all the equalities and inequalities we find that indeed $x^{T}S_{e}^{+}x\leq x^{T}\hat{S}_{e}^{+}x$.
\end{proof}
Recall Theorem~\ref{thm:approx2}, which shows that upper bounds
on the leverages can be used to sample elements and still get an high
quality preconditioner as long as the sample size is increased (in
an easy to compute manner). The last theorem implies that we can find
such an upper bounds using only some of the elements. The crucial
observation is the following corollary to Theorem~\ref{thm:rayleigh}.
\begin{cor}
\label{cor:local}Consider the same conditions as in Theorem~\ref{thm:rayleigh}.
Let $\tilde{\tau}_{e}=\lambda_{\max}(\tilde{K}_{e},\hat{S}_{e})$.
Then we have $\tilde{\tau}_{e}\geq\tau_{e}$. \end{cor}
\begin{proof}
Follows from the previous theorem and the fact that $\lambda_{\max}(\tilde{K}_{e},S_{e})=\lambda_{\max}(S_{e}^{+},\tilde{K}_{e}^{+})$
and $\lambda_{\max}(\tilde{K}_{e},\hat{S}_{e})=\lambda_{\max}(\hat{S}_{e}^{+},\tilde{K}_{e}^{+})$.
\end{proof}
Consider a subset $\hat{E}\subseteq[m]$ of the elements, and let
\[
\hat{K}_{e}=\begin{cases}
K_{e} & e\in\hat{E}\\
\alpha K_{e} & e\notin\hat{E}
\end{cases}
\]
for some $\alpha\in(0,1]$. The last corollary asserts that $\hat{\tau}_{e}\geq\tau_{e}$.
Let $L$ be equal to $\sum_{e\in\hat{E}}K_{e}$ restricted to to non-zero
indexes ($\cup_{e\in\hat{E}}{\cal N}_{e}$). As long as $L$ is well-formed
as well, taking $\alpha\rightarrow0$ and using a continuity argument
we find that the leverage of $e\in\hat{E}$ inside $L$ is an upper
bound to the leverage of $e$ in $K$. However, $L$ might contain
much less elements, so computing the leverage of $e$ inside it might
be cheaper.

This suggest the following local approximation scheme: for an element
$e$, use the effective-stiffness formulas on an element $e$ and
the elements within some distance from it (instead of the entire finite-element
mesh). As argued, this yields an upper bound $\tilde{\tau}_{e}\geq\tau_{e}$.
For this bound to be useful we also need it not to be too loose (otherwise
a huge number of elements will have to be sampled). While we are unable
to characterize exactly when the bound will be loose, and when not,
intuitively a loose bound for an element corresponds to many global
(as opposed to local) behaviors affecting an element, and there are
not too many such global behaviors in a typical finite element models
from applications. Notice that only a small number of loose upper
bounds will not be too detrimental when applying Theorem~\ref{thm:rayleigh}.
Another issue is that we need to compute the leverages for every element.
For this too be cheap we need small local matrices. Again, for finite
element applications who typically have not too-complex geometry this
will typically be the case. Therefore, we believe that this is is
an effective method for sparsifying large meshes.

\section{Numerical Experiments }

In this section we describe two small numerical experiments. Our goal
is to explore how the leverages look on actual finite element matrices,
and show that effective stiffness sampling can indeed select a subset
of the elements to obtain an high-quality preconditioner. We do not
claim to present full practical solver. As we explain in the next
section, and see in the experiment, there are a few challenges that
need to be addressed for that.

\begin{figure}
\noindent \begin{centering}
\begin{tabular}{ccc}
\includegraphics[width=0.45\textwidth]{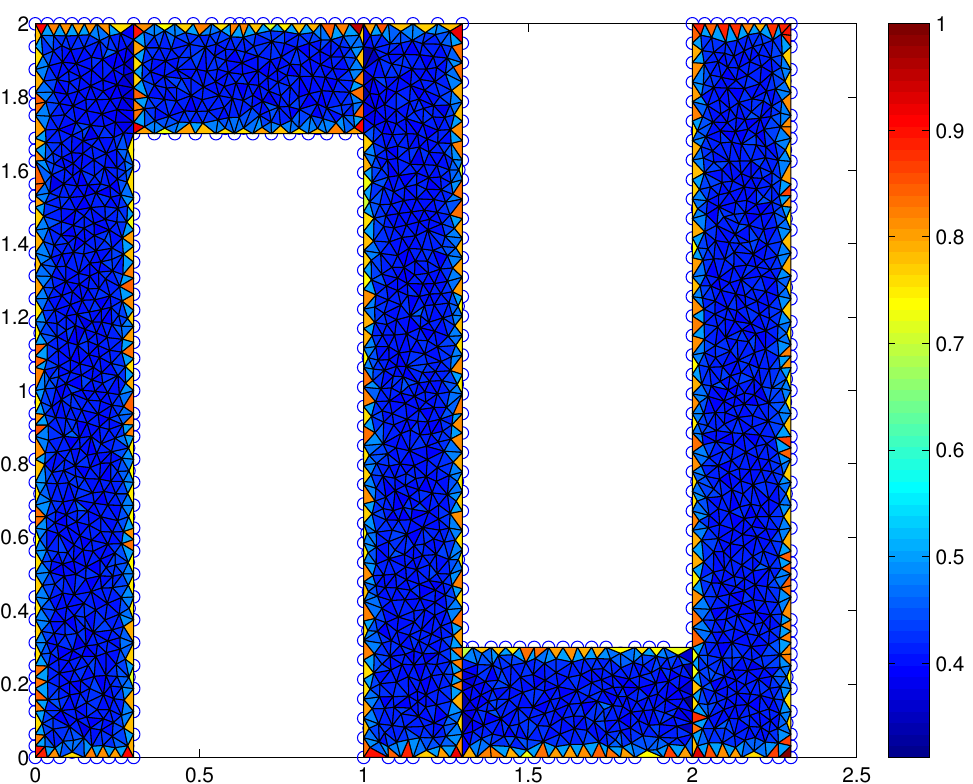} & ~ & \includegraphics[width=0.45\textwidth]{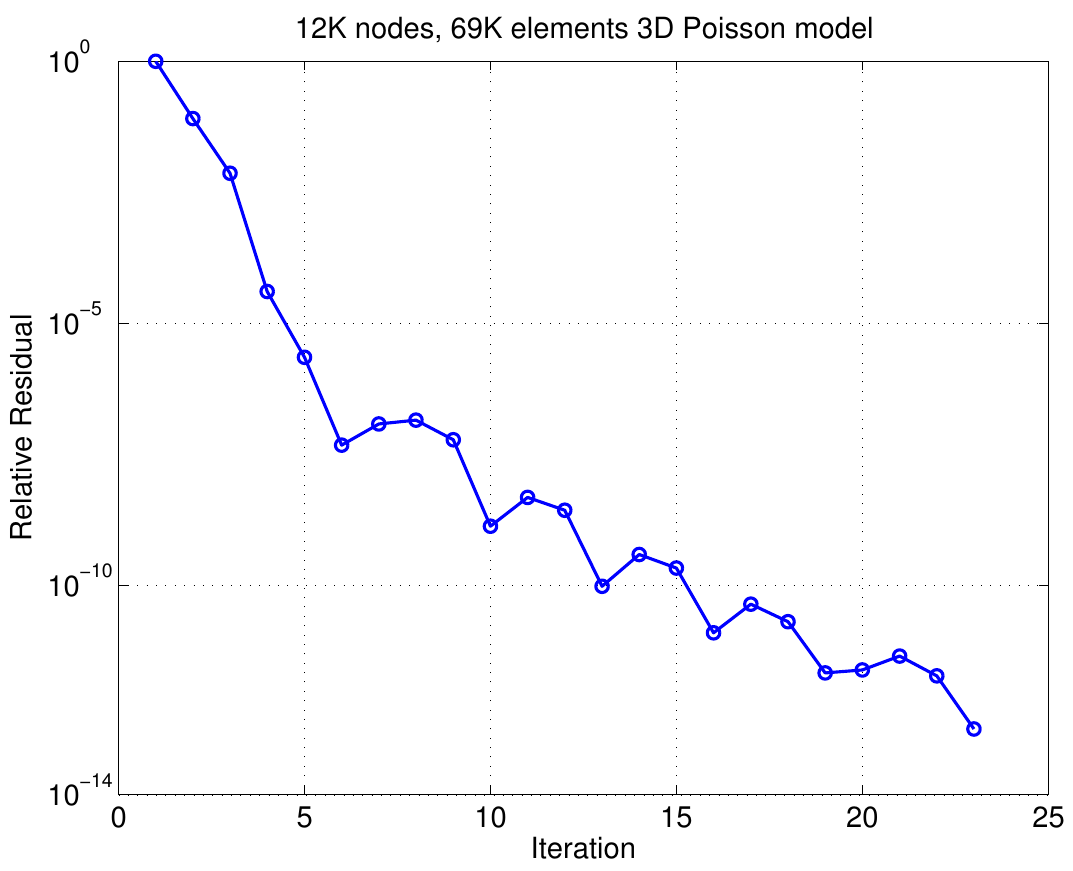}\tabularnewline
\end{tabular}
\par\end{centering}

\caption{\label{fig:example}A numerical example of finite element sparsification.
In the left graph we see S-shaped domain discretized using a triangulated
mesh. The color of each element codes the approximate leverages computed
using a small radius around the element. The right graph shows the
residual as a function of the CG iteration number where the preconditioner
is formed based on effective stiffness sampling. The number of elements
sampled is $ct\log(t)$ for $c=1,2,3,4$ where $t$ is the sum of
approximate leverages.}
\end{figure}

In the first experiment we consider a 2D linear-elasticity problem
on a S-shaped domain discretized using a triangulated mesh. See the
left side of Figure~\ref{fig:example}. There are 1898 nodes, and
3487 elements. The essential element matrices are of size $6$-by-$6$.
The two horizontal bars have significantly different material coefficient
then the three verticals bars (one much weaker, and one much stronger). 

To approximate the leverages we used the local approximation described
in the previous section. For every element we found all the elements
at distance at most 2 from it in the rigidity graph (see~\cite{ShklarskiToledo08}).
Using the rigidity graph ensures we are getting a rigid sub-model.
We computed the effective stiffness matrix of the element inside that
sub-model, and used the approximate stiffness matrix to compute approximate
leverages. The average number of nodes in the local sub-models is
$22$, and the maximum is $24$, so the cost of approximating the
leverage score of an element is about the same as the cost of factoring
a $22$-by-$22$ matrix. Corollary~\ref{cor:local} ensures we are
getting a an upper bound on the leverages. The left side of Figure~\ref{fig:example}
color codes the different elements according to the leverages. We
see that the approximate leverages indeed capture (by giving high
leverages) the important parts of the model: the outer boundary (which
is critical) and the interface between different materials. 

The sum of approximate leverages $\tilde{\tau}_{K}$ is about $1887.9$,
which is about half of the number of nodes (which is the only upper
bound we have on $\tau_{K}$). $\tilde{\tau}_{K}$ can be shrank by
using a smaller radius, e.g. when using a radius $5$ for computing
the approximate leverage scores, the sum $\tilde{\tau}_{K}$ drops
to about $1605.0$. Less elements need to be sampled with this value.
However, the average number of nodes in the local sub-models will
increase to about 63, and the maximum to 84. The time it takes to
approximate the leverages will increase accordingly, so there is clear
a trade-off here.

Theorem~6.2 relates $\tilde{\tau}_{K}$ to the required sample size.
If we apply it to this problem, even when sampling exactly $\left\lceil \tilde{\tau}_{K}\log(\tilde{\tau}_{K})\right\rceil $,
which is below the required number, we sample nearly all of the elements.
It is then no wonder that we get a very good preconditioner. 

We therefore explore convergence in a second experiment. We consider
a synthetic 3D Poisson model with linear elements (essential element
matrices are $4$-by-$4$). The model consists of a ball of one material
inside a box of another material. The model has 12,367 nodes and 69,405
elements. We again compute approximate leverage scores using radius
2 local matrices (average size of local matrices is about $160$-by-$160$).
The approximate leverage sum $\tilde{\tau}_{K}$ is about 2/3 of the
number of nodes.

We now tested convergence of CG when the preconditioner is obtained
using a sample size of $\left\lceil \tilde{\tau}_{K}\log(\tilde{\tau}_{K})\right\rceil $.
We found that convergence is very fast, between 15 to 30 iterations
in all our runs. See the right graph in Figure~\ref{fig:example}
for a typical behavior of the residual. This indicates that the condition
number is not too large. Interestingly, we are sampling \emph{less
}then what is required by the theorems, so the bound seems to be rather
loose. However, when sampling less than $\left\lceil \tilde{\tau}_{K}\log(\tilde{\tau}_{K})\right\rceil $
we frequently got rank-deficient preconditioners.

The value of $\left\lceil \tilde{\tau}_{K}\log(\tilde{\tau}_{K})\right\rceil $
is actually larger than the number of elements in the model. However,
the probabilities are skewed, and the sampling is done with replacement,
so some elements are sampled again and again. It turns out that only
about 50\% of the elements appear in the sampled model. 

We also tried to sample $\left\lceil \tilde{\tau}_{K}\log(\tilde{\tau}_{K})\right\rceil $
elements using uniform samples. Despite the fact that we end up with
more elements (about 65\% of the elements are kept), the sampled model
always lost rank compared to the original one (so it cannot be used
as preconditioner). It seems that non-uniform sampling is essential,
and that the leverage scores provide the necessary probabilities.

\section{\label{sec:discussion}Discussion and Conclusions}

The results in this paper do not constitute practical solver. What
are the remaining challenges that need to be addressed to construct
a complete solver?
\begin{itemize}
\item \textbf{Computing the leverages.} Neither of the two formulas for
the leverages can be computed more efficiently than solving the linear
system itself. Theorems~\ref{thm:approx1} and \ref{thm:approx2}
show that an approximation or an upper bound of the true leverages
suffice. We described a local approximation scheme that might be effective
for large meshes, but further analysis is necessary.
\item \textbf{Number of elements in the sparsified system. }Theoretically,
the number of elements in an $n$-by-$n$ finite element matrix can
be as large as $\Theta(n^{d})$, in which case $O(n\log n)$ elements
is a big improvement. In practice, there are typically only $O(n)$
elements, so sampling $O(n\log n)$ element is not an improvement.
It is worth noting that elements are sampled \emph{with repetition
}so in practice fewer than $O(n\log n)$ distinct elements are sampled.
If the probabilities (leverages) are highly skewed then the number
of element can sampled can even approach $n$. An illustrative, but
unrealistic, example is the following. Consider a finite element matrix
with exactly $n$ elements with leverage $1$ and all other elements
with leverage $0$. All samples will be inside the group of elements
with leverage $1$, so there will only be $n$ distinct elements in
the sample. The authors of \cite{KoutisMillerPang10,KoutisMillerPang11}
used highly skewed probabilities to handle SDD matrices with only
$O(n)$ non-zeros.
\item \textbf{Non-zeros in factor. }Once the finite element matrix has been
sparsified, the sparsified matrix has to be factored, or it can serve
as a foundation for a multilevel scheme. The cost of factoring sparsified
the matrix and the cost of each iteration of PCG depend mainly on
the number of non-zeros in the factor (the fill-in) and not on the
number of non-zeros in the sparsified matrix. To build an effective
preconditioner using sampling, the sampling must be guided so that
the sampled matrix will have low fill. For the sparsifier to be useful
in a multilevel scheme, the sparsified matrix must be easy to coarsen
(eliminate vertices, faces, or elements to obtain a small mesh on
which the process can be repeated).
\end{itemize}
The same issues prevented the initial theoretical results of~\cite{SpielmanSrivastava08}
from immediately producing a fast algorithm. But a few years later
fast algorithms based of effective resistance sampling were suggested.
The extension of effective resistance to effective stiffness is not
trivial. We should not expect the other techniques used in SDD solvers
to trivially extend to finite-element matrices. For example, the first
step in the fastest known SDD solver~\cite{KoutisMillerPang11} is
forming a low-stretch tree. There is currently no equivalent combinatorial
object for finite-element matrices.

Hopefully, our first step will be followed by additional ones that
will enable the construction of general and efficient finite-element
solvers.

\section*{Acknowledgments}

Haim Avron acknowledges the support from XDATA program of the Defense
Advanced Research Projects Agency (DARPA), administered through Air
Force Research Laboratory contract FA8750-12-C-0323. Sivan Toledo
was supported by grant 1045/09 from the Israel Science Foundation
(founded by the Israel Academy of Sciences and Humanities) and by
grant 2010231 from the US-Israel Binational Science Foundation.

\bibliographystyle{plain}
\bibliography{resist-fem}

\end{document}